\documentclass[11pt]{amsart}

\usepackage{bm}
\usepackage{fullpage}
\usepackage{amssymb}
\usepackage{amsmath,amsfonts,amsthm}
\usepackage{hyperref}
\usepackage{color}
\usepackage{cite}
\usepackage{tikz}

\newtheorem{theorem}{Theorem}[section]

\newtheorem{observation}{Observation}[theorem]

\newtheorem{lemma}[theorem]{Lemma}

\newtheorem{claim}[theorem]{Claim}
\newtheorem{proposition}[theorem]{Proposition}

\title{The size, multipartite Ramsey numbers for $nK_2$ versus path-path and cycle }
\author[1]{Yaser Rowshan$^1$}
\author[2]{Mostafa Gholami$^1$}

\address{$^1$Department of Mathematics, Institute for Advanced Studies in Basic Sciences (IASBS), Zanjan 66731-45137, Iran}
\email{y.rowshan@iasbs.ac.ir}
\email{gholami.m@iasbs.ac.ir}

\begin{document}
\maketitle

\begin{abstract}
 For given graphs $G_1, G_2,\ldots, G_n$ and any integer $j$, the size of the multipartite Ramsey number $m_j(G_1, G_2,\ldots, G_n)$ is the smallest positive integer $t$ such that any $n$-coloring of the edges of  $K_{j\times t}$ contains a monochromatic copy of $G_i$ in color $i$ for some $i$, $1 \leq i \leq n$, where $K_{j\times t}$ denotes the complete multipartite graph having $j$ classes with $t$ vertices per each class. In this paper we compute the size of the multipartite  Ramsey number $m_j(K_{1,2}, P_4, nK_2)$ for any $j,n\geq 2$ and  $m_j(nK_2,C_7)$, for any $j\leq4$ and $n\geq 2$.
\end{abstract}

\section{Introduction}
\label{intro}
In this paper, we are only concerned with undirected, simple and finite graphs. We follow \cite{bondy1976graph} for terminology and notations not defined here. For a given graph $G$, we denote its vertex set, edge set, maximum degree and minimum degree by $V(G)$, $E(G)$, $\Delta(G)$ and $\delta(G)$, respectively. For a vertex $v\in V(G)$, we use $\deg_G{(v)}$  and $N_G(v)$  to denote the degree and neighbours of $v$ in $G$, respectively. The neighbourhood of a vertex $v\in V(G)$ are denoted by $N_G(v)=\lbrace u \in V(G) ~ \vert ~ uv \in E(G) \rbrace$ and $N_{X_j}(v)=\lbrace u \in V(X_j) ~\vert~ uv \in E(G) \rbrace$. \\ As usual, a cycle and a path on $n$ vertices are denoted by $C_n$ and $P_n$, respectively. Also by a stripe $mK_2$ we mean a graph on $2m$ vertices and $m$ independent edges. As usual the complete multipartite graph with the partite set $(X_1,X_2,\ldots X_j)$, $|X_i|=s$ for $i=1,2,\ldots j$, denoted by $K_{j\times s}$. We use $[X_i,X_j]$ to denote the set of edges between partite sets $X_i$ and $X_j$. The complement of a graph $G$, denoted by ${\overline G}$, is a graph with the same vertices as $G$ and contains those edges which are not in $G$. Let $T\subseteq V(G)$ be any subset of vertices of G. Then the induced subgraph G[T] is the graph whose vertex set is T and whose edge set consists of all of the edges in E(G) that have both endpoints in T.

\bigskip
For given graphs $G_1, G_2,\ldots, G_n$ and integer $j$, the size of the multipartite Ramsey number $m_j(G_1, G_2,\ldots, G_n)$ is the smallest integer $t$ such that any $n$-coloring of the edges of  $K_{j\times t}$ contains a monochromatic copy of $G_i$ in color $i$ for some $i$, $1 \leq i \leq n$, where $K_{j\times t}$ denotes the complete multipartite graph having $j$ classes with $t$ vertices per each class. $G$ is $n$-colorable to $(G_1, G_2,\ldots, G_n)$ if there exist a $t$-edgs decomposition of $G$ say $(H_1, H_2,\ldots, H_n)$ where $G_i\nsubseteq H_i$ for each $i=1,2, \ldots,n.$\\

 The existence of such a positive integer is guaranteed by a result in \cite{erdos2009partition}. The size of the multipartite Ramsey numbers of small paths versus certain classes of graphs have been studied  in \cite{luczak2018multipartite,sy2011size,jayawardene2016size}. The size of the multipartite Ramsey numbers of stars versus certain classes of graphs have been studied  in  \cite{yek,lusiani2017size}.
In \cite{burger2004ramsey,burger2004diagonal}, Burger,  Stipp, Vuuren, and Grobler, investigated the multipartite ramsey numbers $m_j(G_1, G_2)$, where $G_1$ and $G_2$ are complete balanced multipartite graph, which can be naturally extended to several colors. Recently the numbers $m_j(G_1, G_2)$ have been investigated for special classes: stripes versus cycles, stars versus cycles, see \cite{jayawardene2016size} and its references.\\ In \cite{lusiani2017size}, Lusiani et al. determined the size multipartite Ramsey numbers of $m_j(K_{1,m},H)$, for $j = 2,3$ where $H$ is a path or a cycle on $n$ vertices, and $K_{1,m}$ is a star of order $m + 1$. In this paper we compute the size of the multipartite  Ramsey number $m_j(K_{1,2}, P_4, nK_2)$ for $n,j\geq 2$ and  $m_j(nK_2,C_7)$, for $j\leq4$ and $n\geq 2$ wich are new results of multipartite ramsey numbers. In particular, as the frst target of this work we prove the following theorem:
\begin{theorem}\label{th1}
	$m_j(K_{1,2}, P_4, nK_2)=\lfloor \frac{2n}{j} \rfloor+1$ where $j,n\geq 2$.
\end{theorem}
 In \cite{jayawardene2016size}, Jayawardene, et al. determined the size multipartite Ramsey numbers $m_j(nK_2, C_m)$ where $j\geq 2$ and $m\in \{3,4,5,6\}$. The second goal of this work generalizes these results, as stated below. 
\begin{theorem}\label{th2}
	Let $j\in\{ 2,3,4\}$ and $n\geq 2$. Then 
	\[
	m_j(nK_2,C_7)= \left\lbrace
	\begin{array}{ll}
	\infty & ~~~~~~j=2,n\geq2 ,~\vspace{.2 cm}\\
	2 & ~~~~~~(j,n)=(4,2),~\vspace{.2 cm}\\
	3 & ~~~~~~(j,n)=(3,2), (4,3), ~\vspace{.2 cm}\\
	n & ~~~~~~j=3,n\geq 3 ,~\vspace{.2 cm}\\
	\lceil \frac{n+1}{2} \rceil &  ~~~~~~ j=4,n\geq 4 .
	\end{array}
	\right.
	\]
\end{theorem}

\section{Proof of Theorem \ref{th1}}
In order to simplify the comprehension, let us split the proof of Theorem \ref{th1} into  small parts. We begin with a simple but very useful general lower bound in the following theorem:

\begin{theorem}\label{th3}
	$m_j(K_{1,2}, P_4, nK_2)\geq\lfloor \frac{2n}{j} \rfloor+1$ where $j, n\geq 2$.
\end{theorem}
\begin{proof}
	Consider $G=K_{j\times t}$ where $t=\lfloor \frac{2n}{j} \rfloor$  with partition sets $X_i$, $X_i=\{x^i_1,x^i_2,\ldots ,x^i_t\}$ for $i\in \{1,2,\ldots,j\}$. Consider $x^1_1\in X_1$, decompose the edges of $K_{j\times t}$ into graphs $G_1,G_2$, and $G_3$, where $ G_{1} $ is a null graph and $G_2=\overline{G_3}$, where $ G_{3} $ is  $G[X_1\setminus\{ x^1_1\},X_2,\ldots,X_j]$. In fact $G_{2}$ is isomorphic to $K_{1,(j-1)t}$ and:
	\[ E(G_2) = \{ x^1_1x_i^r~|~~r=2,3,\ldots,j ~\text{and}~ i=1,2\ldots, t\}\] 
	Clearly $E(G_t) \cap E(G_{t'})=\emptyset$, $E(G)=E(G_1) \cup E(G_2)\cup E(G_3)$, $ K_{1,2} \not\subseteq G_1$ and $ P_4 \not\subseteq G_2$. Since $|V(K_{j\times t})|=j\times \lfloor \frac{2n}{j} \rfloor \leq 2n$, we have $|V(G_3)|\leq 2n-1$, that is, $ nK_{2} \not\subseteq G_3$, which means that $m_3(K_{1,2}, P_4, nK_2)\geq \lfloor \frac{2n}{j} \rfloor+1$ and the proof is complete. 
\end{proof} 

\begin{observation}
	\label{ob1} 
	Let $G= K_{2,3}  ( ~or~ K_4-e)$. For any subgraph of $G$, say $H$, either $H$ has a subgraph isomorphism to $K_{1,2}$ or $\overline{H}$ has a subgraph isomorphism to $P_4$.
\end{observation}
\begin{proof}
	Let $H\subseteq G= K_{2,3}$, for $G= K_4-e$ the proof is same. Without loss of generality (W.l.g) let $X=\{x_1,x_2\}$ and $Y=\{y_1,y_2,y_3\}$ be a partition set of $V(G)$ and $P$ be a maximum path in $H$. If $|P|\geq 3$, then $H$ has a subgraph isomorphic to $K_{1,2}$, so let $|P|\leq 2$. If $|P|=1$, then $\overline{H}  (=G)$ has a subgraph isomorphic to $P_4$. Hence we may assume that $|P|=2$, w.l.g let $P=x_1y_1$. Since $|P|=2$, $x_1y_2, x_1y_3$ and $x_2y_1$ are in $E(\overline{H})$ and there is at least one edge of $\{x_2y_2, x_2y_3\}$ in $\overline{H}$, in any case, $P_4\subseteq \overline{H}$ and the proof is complete.
\end{proof}
  
 \bigskip
  We determined the exact value of the multipartite Ramsey number of $m_2(K_{1,2}, P_4, nK_2)$ for $n\geq 2$ in the following lemma:
 
 \begin{lemma}\label{l1}
 	$m_2(K_{1,2}, P_4, nK_2)= n+1$ for $n\geq 2$.
 \end{lemma}
 \begin{proof}
 	Let $X=\{x_1,x_2,\ldots,x_{n+1}\}$ and $Y=\{y_1,y_2,\ldots,
 	y_{n+1}\}$ be a  partition set of $G=K_{n+1,n+1}$. Consider a 3-edge coloring $G^{r}$, $G^{b}$ and $G^{g}$ of $G$. By Theorem \ref{th3} the lower bound holds. Now  let $M$ be the maximum matching in $G^g$. If $|M|\geq n$, then the lemma holds, so let $|M|\leq n-1$. If $|M|\leq n-2$, then we have $K_{3,3}\subseteq \overline{G^g}$ and by Observation \ref{ob1} the lemma holds, so let $|M|=n-1$. W.l.g we may assume that $M = \{x_1y_1, x_2y_2,\ldots,x_{n-1}y_{n-1}\}$. By considering the edges between $\{x_n,x_{n+1}\}$ and $Y\setminus \{y_n,y_{n+1}\}$ and the edges between $\{y_n,y_{n+1}\}$ and $X\setminus \{x_n,x_{n+1}\}$, we have $K_{3,2}\subseteq G^r\cup G^b$. Hence by Observation \ref{ob1} the lemma  holds.
 \end{proof}
 In the next two lemmas, we consider $m_3(K_{1,2}, P_4, nK_2)$ for certain values of $n$. In particular, we prove that 
 $m_3(K_{1,2}, P_4, nK_2)= n$ for $n=2,3$ in Lemma \ref{vlem1} and $m_3(K_{1,2}, P_4, 4K_2)=3$  in Lemma \ref{vlem2}.
 \begin{lemma}\label{vlem1}
 	$m_3(K_{1,2}, P_4, nK_2)= n$ for $n=2,3$.
 \end{lemma}
 \begin{proof}
 	Let $X_i=\{x^i_1,x^i_2,\ldots ,x^i_n\}$ for $i\in \{1,2,3\}$ be a  partition set of $G=K_{3\times n}$. Consider a 3-edge coloring $G^{r}$, $G^{b}$ and $G^{g}$ of $G$. By Theorem \ref{th3} the lower bound  holds. Now  let $M$ be the maximum matching in $G^g$ and consider the following cases:
 	
 	\bigskip
 	{\bf Case 1:} $n=2$. If $|M|\geq 2$ then $nK_2\subseteq G^g$ and the proof is complete. So let $|E(M)|\leq 1$. W.l.g we may assume that $x_1^1x_1^2 \in E(M)$,  hence we have $K_4-e\cong G[x_2^1,x^2_2,X_3]\subseteq G^r\cup G^b$ and by Observation \ref{ob1} the proof is compelet.
 	
 	\bigskip
 	{\bf Case 2:} $n=3$. In this case, if $|E(M)|\leq 1$ or $|E(M)|\geq 3$, then the proof is same as case 1. So let $|E(M)|= 2$ and w.l.g we may assume that $E(M)=\{e_1, e_2\}$. Considering any $e_1$ and $e_2$ in $E(G)$. In any case we have $G^r\cup G^b$  has a subgraph isomorphic to $K_{3,2}$, hence by Observation \ref{ob1} the lemma holds. Therefore we have $m_3(K_{1,2}, P_4, 3K_2)= 3$.\\
 	Now by cases 1 and 2, the proof is complete.
 \end{proof}
 
 \begin{lemma}\label{vlem2}
 	$m_3(K_{1,2}, P_4, 4K_2)=3$.
 \end{lemma}
 \begin{proof}
 	Let  $X_i=\{x^i_1,x^i_2, x^i_3\}$ for $i\in \{1,2,3\}$ be a  partition set of $G=K_{3\times 3}$. By Theorem \ref{th3} the lower bound holds. Consider a 3-edge coloring ($G^{r}$, $G^{b}$, $G^{g}$) of $G$ where $4K_2\nsubseteq G^g$.  Let $M$ be a maximum matching in $G^g$, if $|M|\leq 2$ then the proof is same as Lemma \ref{vlem1}. Hence we may assume that $|M|=3$ and w.l.g let $E(M) = \{e_1,e_2,e_3\}$. By Observation \ref{ob1} there is at least one edge between $X_1$ and $X_2$ in $G^g$, say $e_1=x_1^1x_1^2$, similarly  there is at least one edge between $X_3$ and $\{x_2^1, x_3^1\}$ in $G^g$, say $e_2=x_2^1x_1^3$, otherwise $K_{3,2}\subseteq G^r\cup G^b$ and  the proof is complete. Now by Observation \ref{ob1} there is at least one edge between $\{x_3^1,x_2^3,x_3^3\}$ and $\{x_2^2, x_3^2\}$ in $G^g$, let $e_{3}$ be this edge. If $x_3^1\notin V(e_3)$ (say $e_3=x_2^2x_2^3$), then $K_3\subseteq G^r\cup G^b[x_3^1,x_3^2,x_3^3]$.\\
 	Now consider the vertex $x_1^1$ and $x_1^2$, since $|M|=3$ and  $e_1=x_1^1x_1^2$, it is easy to check that $x_1^1 x_3^3, x_1^2 x_3^3 \in E(G^g)$ and $x_1^1 x_3^2, x_1^2 x_3^1 \in E(\overline {G^g})$, otherwise $K_{4} - e \subseteq \overline{G^g}$ and the
 	proof is complete. Similarly we have $x_2^1 x_3^2, x_1^3 x_3^2 \in E(G^g)$ and $x_2^1 x_3^3, x_1^3 x_3^1 \in E(\overline {G^g})$. Now by considering the edges of $G[X_{1},x_1^2,x_3^2,x_1^3,x_3^3]$ it is easy to check that $K_{4} - e \subseteq G^r\cup G^b $ and the lemma holds. Hence we have $x_3^1\in V(e_3)$ (say $e_3=x_3^1x_2^2$), in this case we have  $K_{2,2}\cong G[x_2^2,x_3^2,x_2^3,x_3^3]\subseteq G^r\cup G^b$, otherwise  if there exist at least one edge between $\{x_2^3,x_3^3\}$ and $\{x_2^2, x_3^2\}$ in $G^g$, say $e$, set  $e=e_3$ and the proof is same. Hence by considering the vertex $x_1^1$ and $x_1^2$, since $|M|=3$ and  $e_1=x_1^1x_1^2$, it easy to check that $K_{3,2}\subseteq G^r\cup G^b$ and by Observation \ref{ob1} the proof is complete.  
 \end{proof}
 \begin{theorem}\label{th4}
 	$m_3(K_{1,2}, P_4, nK_2)\leq \lfloor \frac{2n}{3} \rfloor+1$ for each $n\geq 2$.
 \end{theorem}
 \begin{proof}
 	Let  $X_i=\{x^i_1,x^i_2,\ldots ,x^i_t\}$ for $i\in \{1,2,3\}$ be a partition sets of $G=K_{3\times t}$ where $t=\lfloor \frac{2n}{3} \rfloor+1$. We will prove this theorem by induction. For the base step of the induction, since $\lfloor \frac{2\times 2}{3} \rfloor+1=2, \lfloor \frac{2\times 3}{3} \rfloor+1=3$ and $\lfloor \frac{2\times 4}{3} \rfloor+1=3$, theorem holds by Lemmas \ref{vlem1} and \ref{vlem2}. Suppose that $n\geq 5$ and $m_3(K_{1,2}, P_4, n'K_2)\leq \lfloor \frac{2n'}{3} \rfloor+1$ for each $n'<n$. We will show that $m_3(K_{1,2}, P_4, nK_2)\leq \lfloor \frac{2n}{3} \rfloor+1$. By contradiction, we may assume that $m_3(K_{1,2}, P_4, nK_2)> \lfloor \frac{2n}{3} \rfloor+1$, that is, $K_{3\times(\lfloor \frac{2n}{3}\rfloor+1)}$ is 3-colorable to $(K_{1,2}, P_4, nK_2) $. Consider a 3-edge coloring $(G^{r}, G^{b}, G^g)$ of $G$, such that $K_{1,2} \not\subseteq G^{r}$, $P_4 \not\subseteq G^{b}$ and $nK_2 \not\subseteq G^{g}$. By the induction hypothesis and Theorem \ref{th3} we have $m_3(K_{1,2}, P_4, (n-1)K_2)= \lfloor \frac{2(n-1)}{3} \rfloor+1\leq \lfloor \frac{2n}{3} \rfloor+1 $. Therefore since $K_{1,2} \not\subseteq G^{r}$ and $P_4 \not\subseteq G^{b}$ we have $(n-1)K_2 \subseteq G^g$. Now we have the following cases:
 	
 	\bigskip
 	{\bf Case 1:} $\lfloor \frac{2n}{3} \rfloor = \lfloor \frac{2(n-1)}{3} \rfloor+1$.\\
 	Since $\lfloor \frac{2n}{3} \rfloor = \lfloor \frac{2(n-1)}{3} \rfloor+1$
 	we have a copy of $H=K_{3\times (\lfloor \frac{2(n-1)}{3}\rfloor+1)}$ in $G$. In other words, for each $i\in\{1,2,3\}$, there is a vertex, say $x\in X_i$, such that $x \in V(G)\setminus V(H)$.  W.l.g we may assume that $A=\{x_1^1,x_1^2,x_1^3\}$ be this vertices. Since $H\subseteq G$, we have $K_{1,2} \not\subseteq G^{r}[V(H)]$ and $P_4 \not\subseteq G^{b}[V(H)]$. Hence by the induction hypothesis, we have $M=(n-1)K_2\subseteq G^g[V(H)]\subseteq G^g$. We consider that the three vertices do not belong to $V (H)$, i.e., $A$. Since $nK_2 \not\subseteq G^g$, we have $G[A]\subseteq G^r\cup G^b$. Now we consider the following Claim:
 	 \begin{claim}\label{c1}
 		$n\in B\cup D$ where $B=\{3k~| ~~k=1,2,...\}$ and $D=\{3k+2~| ~~k=1,2,...\}$.
 	\end{claim}
 	{\it Proof.} By contradiction we may assume that $n\notin B\cup D$. In other words, let $n=3k+1$, then we have
 		\[2k=\lfloor \frac{6k}{3} \rfloor=\lfloor \frac{6k}{3}+\frac{2}{3} \rfloor=\lfloor \frac{6k+2}{3} \rfloor=\lfloor \frac{2(3k+1)}{3} \rfloor\]
 		\[=\lfloor \frac{2n}{3} \rfloor=\lfloor \frac{2(n-1)}{3} \rfloor+1=\lfloor \frac{2(3k)}{3} \rfloor+1=2k+1,\]
 		which is  a contradiction implying that $n\in B\cup D$.	
 	
 	\begin{claim}\label{c2}
 		There is at least one vertex in $ V(H)\setminus V(M)$.
 	\end{claim}
{\it Proof.} Let $M=(n-1)K_2\subseteq G^g$, Then $|V(M)|=2(n-1)=2n-2$. Since $\lfloor \frac{2n}{3} \rfloor = \lfloor \frac{2(n-1)}{3} \rfloor+1$, by Claim 1, if $n\in B$ we have $n=3k$ for $ k\geq 2$. Now we have
 		\[ \lfloor \frac{2(n-1)}{3} \rfloor+1=\lfloor \frac{2(3k-1)}{3} \rfloor+1=\lfloor \frac{2(3k)}{3}-\frac{2}{3} \rfloor+1=2k-1+1=2k.\]
 		Hence we have $|V(H)|=3\times(2k)=6k=2n$ and thus $|V(H)\setminus V(M)|=2$. If $n\in D$  then we have
 		\[ \lfloor \frac{2(n-1)}{3} \rfloor+1=\lfloor \frac{2(3k+1)}{3} \rfloor+1=\lfloor \frac{2(3k)}{3}+ \frac{2}{3} \rfloor+1= 2k+1.\]
 		Hence $|V(H)|=3\times(2k+1)=6k+3=2n-1$. Therefore $|V(H)\setminus V(M)|=1$.
 	
 \bigskip	
 	Now by Claim \ref{c2}, let $x\in V(H)\setminus V(M)$. Since $nK_2\not\subseteq G^g$ we have $K_4-e\cong G[A\cup \{x\}]\subseteq G^r\cup G^b$. Hence by Observation \ref{ob1} again a contradiction.
 	
 	\bigskip
 	{\bf Case 2:} $\lfloor \frac{2n}{3} \rfloor = \lfloor \frac{2(n-1)}{3} \rfloor$.\\
 	In this case, by Claim \ref{c1} we have $n=3k+1$. Since $K_{1,2} \not\subseteq G^{r}$ and $P_4 \not\subseteq G^{b}$, by the induction hypothesis we have $M=(n-1)K_2\subseteq G^g$. Now we have the following claim:
 	\begin{claim}\label{c3} 
 		$|V(G)\setminus V(M)|=3$.
 	\end{claim}
 {\it Proof.} Let $M=(n-1)K_2\subseteq G^g$. Since $|V(X_j)|=\lfloor \frac{2n}{3} \rfloor+1$ and $n=3k+1$, we have $\lfloor \frac{2n}{3} \rfloor+1=\lfloor \frac{2(3k+1)}{3} \rfloor+1=\lfloor \frac{6k}{3}+\frac{2}{3} \rfloor+1=2k+1$ and, therefore, $|V(G)|=3\times (2k+1)=6k+3=2(3k+1)+1=2n+1$, that is, $|V(G)\setminus V(M)|=(2n+1)-(2n-2)=3$.
 
 \bigskip
 	By Claim \ref{c3} we have $|V(G)\setminus V(M)|=3$.  W.l.g we may assume that $A'=\{x,y,z\}$ is this vertices, since $nK_2\not\subseteq G^g$ we have $G[A']\subseteq G^r\cup G^b$. We consider the three vertices belonging to $A'$, now we have the following subcases:
 	
 	\bigskip
 	{\bf Subcase 2-1:} $A'\subseteq X_j$ for only one $j\in\{1,2,3\}$.  W.l.g we may assume that $A'\subseteq X_1$ and $E(M)=\{e_i~|~i=1,2,\ldots, (n-1)\}$. Since $k\geq 2$ and $3k+1=n\geq 7$ we have $|X_j|\geq 5$ and $|E(M)\cap E(G[X_2,X_3])|\geq 3$, otherwise, $K_{3,3}\subseteq G^r\cup G^b$ and by Observation \ref{ob1}; a contradiction.  W.l.g we may assume that $\{x_i^2x_i^3 ~|~i=1,2,3\}\subseteq (E(M)\cap E(G^g[X_2,X_3]))$. Consider $G'=G[A',x_1^2,x_2^2,x_3^2,x_1^3,x_2^3,x_3^3]\cong K_{3 \times 3}$. Since $nK_2\not\subseteq G^g$, if $M'$ is a maximum matching in $G'^{g}$, then $|M'|\leq 3$, otherwise we have $nK_2=M\setminus \{e_1,e_2,e_3\}\cup M'\subseteq G^g$; a contradiction again. Since $m_3(K_{1,2},P_4,4K_2)=3$ and $|M'|\leq 3$, we have $K_{1,2} \subseteq G'^r\subseteq G^r$ or $P_4\subseteq G'^b\subseteq G^b$; also a contradiction.
 	
 	\bigskip
 	{\bf Subcase 2-2:} $|A'\cap X_j|=1$ for each $j\in\{1,2,3\}$.  W.l.g we may assume that $x\in X_1, y\in X_2 $ and $z\in X_3$. Hence $G[A']\cong K_3\subseteq G^r\cup G^b$. Since  $|X_j|\geq 5$ we have $|E(M)\cap E(G^g[X_i,X_j])|\geq 2 $ for each $i,j\in \{1,2,3\}$.  W.l.g we may assume that $x'y'\in E(M)\cap E(G^g[X_1\setminus \{x\},X_2 \setminus \{y\}])$, $x'\in X_1$ and $y'\in X_2$. If $ x'y$ and $x'z\in E(G^r\cup G^b)$ then we have $K_4-e\subseteq G^r\cup G^b$ and by Observation \ref{ob1}; a contradiction. So let $ x'y$ or $x'z\in E(G^g)$. If $ x'y\in E(G^g)$, then, since $nK_2\not\subseteq G^g$, we have $y'x, y'z\in E(G^r\cup G^b)$, that is, $K_4-e\subseteq G^r\cup G^b$; a contradiction again. So let $ x'z\in E(G^g)$ and $x'y\in E(G^r\cup G^b)$. Since $nK_2\not\subseteq G^g$, we have $y'x\in E(G^r\cup G^b)$. If $|E(G^r)\cap E(G[A'])|\not=0$ then we have $P_4\subseteq G^b$. So let $xy, yz,zx\in E(G^b)$ and $xy',yx'\in E(G^r)$. Since $|E(M)\cap E(G^g[X_i,X_j])|\geq 2 $ there is at least one edge, say $y''z''\in E(M)\cap E(G^g[X_2\setminus \{y\},X_3\setminus \{z\}])$. W.l.g we may assume that $y''\in X_2$ and $z''\in X_3$. Since $K_{1,2} \not\subseteq G^r$ and $P_4\not\subseteq G^b$ we have $y''x, z''y \in E(G^g)$. Hence we have $nK_2=M\setminus\{y''z''\}\cup \{y''x,z''y\}$; a contradiction.
 	
 	\bigskip
 	{\bf Subcase 2-3:} $|A'\cap X_j|=2$ for only one $j\in\{1,2,3\}$.  W.l.g we may assume that $x, y\in X_1$ and $ z\in X_2$. Hence we have $G'[A']\cong P_3\subseteq G^r\cup G^b$. Since $k\geq 2$ we have $|X_j|\geq 5$, that is, $|E(M)\cap E(G^g[X_2,X_3])|\geq 3 $.  W.l.g we may assume that $vu,v'u'\in E(M)\cap G^g[X_2,X_3]$ where $v,v'\in X_2$ and $ u,u'\in X_3$. Now we have the following claim:
 	
 	\begin{claim}\label{c4}
 		$|N_{G^g}(x)\cap \{v,v'\}|=|N_{G^g}(y)\cap \{v,v'\}|=0$.
 	\end{claim}
{\it Proof.} By contradiction, w.l.g we may assume that $ xv\in E(G^g)$. Since $nK_2\not\subseteq G^g$, we have $yu, zu \in E(G^r\cup G^b)$. Consider $A''=\{y,z,u\}$ and $M'=M\setminus\{vu\}\cup \{xv\}$. Hence  $M'=(n-1)K_2\subseteq G^g$ and $|A''\cap X_j|\not=0$ for each $j\in\{1,2,3\}$; a contradiction to subcase 2-2.\\
 		Now by Claim \ref{c4} we have $K_{2,3}=G[A'\cup\{v,v'\}]\subseteq G^r\cup G^b$. In this case by Observation \ref{ob1} we have $K_{1,2} \subseteq G^r$ or $P_4\subseteq G^b$; a contradiction again.\\
 	Therefore by Cases 1 and 2 we have $m_3(K_{1,2}, P_4, nK_2)\leq \lfloor \frac{2n}{3} \rfloor+1$ for $n\geq 2$. 
 \end{proof}
 Now by Theorems \ref{th3} and \ref{th4}  we have the following theorem:
 \begin{theorem}\label{th5}
 	$m_3(K_{1,2}, P_4, nK_2)= \lfloor \frac{2n}{3} \rfloor+1$ for $n\geq 2$. 
 \end{theorem}
 In the next two theorems, we consider $m_j(K_{1,2}, P_4, nK_2)$ for each values of $n\geq 2$ and $j\geq 4$. In particular, we prove that $m_j(K_{1,2}, P_4, nK_2)= \lfloor \frac{2n}{j}\rfloor+1$ for $n\geq 2$ and $j\geq4$.  We start with the following theorem:
 
 \bigskip
 \begin{theorem}\label{th6}
 	Let $j\geq 4$ and $n\geq 2$. Given that $m_j( K_{1,2}, P_4, (n-1)K_2)= \lfloor \frac{2(n-1)}{j}\rfloor +1$, it follows that $m_j(K_{1,2}, P_4, nK_2) \leq \lfloor \frac{2n}{j}\rfloor+1.$
 \end{theorem}
 \begin{proof} Let $j\geq 4$ and $n\geq 2$. For $i\in \{1,2,\ldots,j\}$ let $X_i=\{x^i_1,x^i_2,\ldots ,x^i_t\}$ be  partition set of $G= K_{j\times t}$ where $t= \lfloor \frac{2n}{j}\rfloor+1$. Assume that $m_j(K_{1,2}, P_4, (n-1)K_2) = \lfloor \frac{2(n-1)}{j}\rfloor +1$ is true. To prove $m_j(K_{1,2}, P_4, nK_2) \leq \lfloor \frac{2n}{j}\rfloor+1$. Consider a 3-edge coloring $(G^{r}, G^{b}, G^{g})$ of $G$. Suppose that $nK_2 \not\subseteq G^{g}$, we prove that $ K_{1,2} \subseteq G^{r}$ or $ P_4 \subseteq G^{b}$. Let $M^*$ be the maximum matching in $G^{g}$. Hence by the assumption, $|M^*|\leq n-1$, that is $|V(K_{j\times t}) \cap V(M^*)|\leq 2(n-1)$. Now, we have the following claim:
 	\begin{claim}\label{c5}
 		$|V(K_{j\times t}) \setminus V(M^*)|\geq 3$.
 	\end{claim}
{\it Proof.} Consider the following cases:
 		
 		\bigskip
 		{\bf Case 1:} Let $2n=jk ~(2n\equiv 0( mod ~ j))$. In this case, we have
 		\[|V(G)|=j\times t=j \times ( \lfloor \frac{2n}{j}\rfloor+1)= j \times\lfloor \frac{2n}{j}\rfloor+j= jk+j=j(k+1).\]
 		Hence
 		\[|V(G)\setminus V(M^*)|\geq j(k+1) -2(n-1)= jk+j-2n+2=j+2 \geq 6    ~ (j \geq 4).\] 
 		
 		\bigskip
 		{\bf Case 2:} Let $2n=jk+r~(2n\equiv r( mod ~ j)$ where $r\in \{1,2,\ldots,j-1\})$. In this case, we have\\ 
 		$|V(G)|=j \times ( \lfloor \frac{2n}{j}\rfloor+1)=j \times ( \lfloor \frac{jk+r}{j}\rfloor+1)=j \times ( \lfloor \frac{jk}{j}+ \frac{r}{j}\rfloor+1)= j \times\lfloor \frac{jk}{j}\rfloor+j= jk+j.$\\
 		Hence we have\\ 
 		$|V(G)\setminus V(M^*)|\geq j(k+1) -2(n-1)= jk+j-2n+2=jk+j-jk-r+2=j-r+2 \geq 3.$
 	
 		\bigskip 	
 	By Claim \ref{c5}, $G$ contains three vertices, say $x,y$ and $z$ in $V(K_{j\times t}) \setminus V(M^*)$.  Consider the vertex set $\{x,y,z\}$ and let $\{x,y,z\}\subseteq A=V(G)\setminus V(M^*)$. Now we have the following cases:
 	
 	\bigskip
 	{\bf Case 1:} Let $x\in X_1$, $y\in X_2$ and $z\in X_3$ where $X_i$ for $i=1,2,3$ are distinct partition sets of $G=K_{j\times t}$. Note that  all vertices of $A$ are adjacent to each other in $\overline{G^g}$. Since $t\geq 2$ we have $|X_i|\geq 2$. Consider the partition $X_j$ for $j\geq 4$. Since $|X_j|\geq 2$, if $|A\cap X_j| \geq 1$ for at least one $j\geq 4$, then we have $K_4\subseteq G^r\cup G^b$ and  the proof is complete by Observation \ref{ob1}. Now let $|A\cap X_j |=0$ for each $j\geq 4$.  Hence for $x^4_1\in X_4$ there exists a vertex, say $u$ such that $x^4_1u\in E(M^*)$. Consider $N_{G^g}(x^4_1)\cap \{x,y,z\}$. If $|N_{G^g}(x^4_1)\cap \{x,y,z\}|\leq 1$, then we have $K_4-e\subseteq G^r\cup G^b$ and, by Observation \ref{ob1}, the proof is complete. Therefore let $|N_{G^g}(x^4_1)\cap \{x,y,z\}|\geq 2$.  W.l.g we may assume that $\{x,y\}\subseteq N_{G^g}(x^4_1)\cap \{x,y,z\}$. In this case, we have $|N_{G^g}(u)\cap \{x,y,z\}|= 0$. On the contrary, let $xu\in E(G^g)$ and set $M'=M^*\setminus \{x^4_1u\}\cup\{x^4_1y,ux\}$.  Clearly $M'$ is a matching where $|M'|>|M^*|$, which contradicts the maximality of $M^*$. Hence we have $|N_{G^g}(u)\cap \{x,y,z\}|=0$. Therefore we have $K_4-e\subseteq G^r\cup G^b[x,y,z,u]$ and,  by Observation \ref{ob1}, the proof is complete. 
 	
 	\bigskip
 	{\bf Case 2:} Let $x,y\in X_i$ and $z\in X_{i'}$ where $X_i, X_{i'}$  are distinct partition sets of $G$.  W.l.g let $i=1$ and $i'=2$.   Consider the partition $X_j ~(j\neq 1,2)$. Since $|X_j|\geq 2$, if $|A\cap X_j| \geq 1$, then we have $K_4-e\subseteq G^r\cup G^b$ and, by Observation \ref{ob1}, the proof is complete. So let $|A\cap X_j |=0$ for each $j\geq 3$.  Now, we have the following claim.
 	
 	\begin{claim}\label{c6}
 		Let $e=v_1v_2\in E(M^*)$, and w.l.g let $|N_{G^g}(v_1)\cap \{x,y,z\}|\geq |N_{G^g}(v_2) \cap \{x,y,z\}|$. If $|N_{G^g}(v_1)\cap \{x,y,z\}|\geq 2$, then $|N_{G^g}(v_2)\cap \{x,y,z\}|=0$. If $|N_{G^g}(v_1)\cap \{x,y,z\}|=|N_{G^g}(v_2)\cap \{x,y,z\}|=1$, then $v_1,v_2$ has the same neighbour in $\{x,y,z\}$.
 	\end{claim}
{\it Proof.} Let $|N_{G^g}(v_1)\cap \{x,y,z\}|\geq 2$. W.l.g we may assume that $\{w,w'\}\subseteq N_{G^g}(v_1)\cap \{x,y,z\}$. By contradiction let $|N_{G^g}(v_2)\cap \{x,y,z\}|\neq 0$, w.l.g  let $w{''}\in N_{G^g}(v_2)\cap \{x,y,z\}$. In this case, we set $ M' = (M^*\setminus \{v_1v_2\}) \cup\{v_1w,v_2w{''}\}$. Clearly $M'$ is a matching with $|M'| > |M^*|$, which contradicts the maximality of $M^*$. So let $|N_{G^g}(v_i)\cap \{x,y,z\}|=1$ for $i=1,2$, if $v_i$ has a different neighbour then the proof is  same.
 	
 	\begin{claim}\label{c7}
 		There is at least one edge, say $e=u_iu_j \in E(M^*)$, such that $u_i, u_j\notin X_1,X_2$.
 	\end{claim}
{\it Proof.} If $|X_j|\geq 3$ then there is at least one edge, say $e=u_iu_j \in E(M^*)$, such that $u_i, u_j\notin X_1,X_2$. Otherwise, we have $K_{3,2}\subseteq G^r\cup G^b[X_j,X_{j'}]$ where $j, j'\geq 3$ , hence by Observation  \ref{ob1}; a contradiction. So let $|X_j|=2$. In this case, if $j\geq 5$ then the proof is same. Now let $j=4$. We have $|M^*|\leq 2$, that is, $n\leq 3$. Hence there is at least one vertex, say $w\in (X_3\cup X_4)\cap A$;  a contradiction to $|A\cap X_j |=0$.
   
 \bigskip  
 	By Claim \ref{c7} there is at least one edge, say $e=u_iu_j \in E(M^*)$, such that $u_i, u_j\notin X_1,X_2$. W.l.g let $e=u_1u_2 \in E(M^*)$ such that $u_i\notin X_1,X_2$, also w.l.g assume that $|N_{G^g}(u_1)\cap \{x,y,z\}| \geq |N_{G^g}(u_2)\cap \{x,y,z\}|$. If $|N_{G^g}(u_1)\cap \{x,y,z\}|\geq 2$, then by Claim \ref{c7} we have $|N_{G^g}(u_2)\cap \{x,y,z\}|=0$. Hence we have $K_4-e\subseteq G^r\cup G^b$. So let $|N_{G^g}(u_1)\cap \{x,y,z\}|=|N_{G^g}(u_2)\cap \{x,y,z\}|=1$, in this case, by Claim \ref{c7}, we have $N_{G^g}(u_1)\cap \{x,y,z\}=N_{G^g}(u_2)\cap \{x,y,z\}$, and if $x$ or $y$ is this vertex, then $K_4-e\subseteq G^r\cup G^b$, otherwise $K_{3,2}\subseteq G^r\cup G^b$. In any case, by Observation \ref{ob1}, the proof is complete.
 	
 	\bigskip
 	{\bf Case 3:} Let $x,y,z\in X_i$ where $X_i$  is a partition set of $G=K_{j\times t}$, say $i=1$. If there exists a vertex, say $w\in X_j \cap A$, where $j\neq1$, then the proof is same as case 2. Hence let $|A\cap X_j|=0$. Since $|X_j|\geq 3$, there exists an edge, say $e=vu \in E(M^*)$, such that $v,u\notin X_1$. Consider the neighbours of vertices $v$ and $u$ in $X_1$. W.l.g let $|N_{G^g}(v)\cap \{x,y,z\}|\geq|N_{G^g}(u)\cap \{x,y,z\}|$. If $|N_{G^g}(v)\cap \{x,y,z\}|=0$ then we have $K_{3,2}\subseteq G^r\cup G^b$, so let $|N_{G^g}(v)\cap \{x,y,z\}|\geq 1$. In this case, by Claim \ref{c7}, we have $|N_{G^g}(u)\cap \{x,y,z\}|\leq 1 $. Hence  w.l.g we may assume that $yu$ and $zu $ be in $E(G^r\cup G^b)$ and $x\in N_{G^g}(v)$. Now set $M^{**}=(M^*\setminus\{vu\})\cup\{vx\}$ and $A'=(A\setminus\{x\})\cup \{u\}$, the proof is same as case 2 and the proof is complete.\\
 	According to the cases $1,2$ and $3$ we have $m_j(K_{1,2}, P_4, nK_2) \leq \lfloor \frac{2n}{j}\rfloor+1$.
 \end{proof}
   The results of Theorems \ref{th3}, \ref{th5} and \ref{th6} and Lemmas \ref{l1}, \ref{vlem1} and \ref{vlem2} concludes the proof of Theorem \ref{th1}.
 %\begin{theorem}\label{vth4}
 %	$m_j(K_{1,2}, P_4, nK_2)= \lfloor \frac{2n}{j} \rfloor+1$, for $j\geq 2$ and $n\geq 2$
 %\end{theorem}
 \section{Proof of Theorem \ref{th2}}
 In this section, we investigate the size multipartite Ramsey number $m_j(nK_2, C_7)$  for $j\leq4$ and $n\geq 2$. In order to simplify the comprehension, let us split the proof of Theorem \ref{th2} into  small parts. For $j=2$, since the bipartite graph has no odd cycle, we have $m_2(nK_2, C_7)=\infty$. For other cases we start with the following proposition:
 \begin{proposition}\label{pr1}
 	$m_3(nK_2,C_7) =3$ where $n=2,3$.
 \end{proposition}
 \begin{proof}
 	Clearly $m_3(nK_2,C_7) \geq3$.  Consider $K_{3\times 3}$ 
 	with the partition set $X_i=\{x_1^i,x_2^i, x_3^i\}$ for $i=1,2,3$. Let $G$ be a subgraph of $K_{3\times3}$. For $n=2$ if $2K_2\subseteq G$, then the proof is complete, so let $2K_2\not\subseteq G$. In this case we have $K_{3,2,2}\subseteq \overline{G}$, hence $C_7\subseteq \overline{G}$, that is, $m_3(2K_2,C_7) =3$. For $n=3$ by contradiction, we may assume that  $m_3(3K_2,C_7) > 3$, that is, $K_{3\times 3}$ is $2$-colorable to $(3K_2,C_7)$, say $3K_2\not\subseteq G$ and $C_7\not\subseteq\overline{G}$. Since $m_3(3K_2,C_6) =3$ \cite{jayawardene2016size}, and $3K_2\not\subseteq G$ we have $C_6\subseteq\overline{G}$. Let $A=V(C_6)$ and $Y_i=A\cap X_i$ for $i=1,2,3$. If there exists $i\in \{1,2,3\}$ such that $|Y_i|=0$, say $i=1$, then we have $A=X_2\cup X_3$ and $C_6\subseteq\overline{G}[X_2,X_3]$. Let $C_6=w_1w_2\ldots w_6w_1$. Since $C_7\not\subseteq\overline{G}$, for each $x_i\in X_1$  in $\overline{G} $, $ x_i$ cannot be adjacent to $w_i$ and $w_{i+1}$ for $i=1,2,\ldots, 6$. Hence we have $|N_G(x_i)\cap V(C_6)|\geq 3$ for each $x_i\in X_1$. One can easily check that in any case,  we have $3K_2\subseteq G$; a contradiction, hence let $|Y_i|\geq 1$ for each $i=1,2,3$. Set $B=(|Y_1|, |Y_2|,|Y_3|)$. Now we have the following cases:
 	
 	\bigskip
 	{\bf Case 1:} $B=(3,2,1)$. let $A=X_1\cup\{x_1^2,x^2_2,x_1^3\}$. In this case, we have $C_6\cong x_1^1x_1^2x_2^1x_2^2x_3^1x_1^3x_1^1$. Consider the vertex set $A' = V(K_{3\times 3})\setminus A=\{x_3^2,x_2^3,x_3^3\}$. Since $C_7\not\subseteq\overline{G}$, we have $|N_{\overline{G}}(x_2^3)\cap \{x_1^1,x_1^2\}|\leq 1$. Hence  $|N_G (x_2^3)\cap \{x_1^1,x_1^2\}|\geq 1$. W.l.g let $x_2^3x_1^1\in E(G)$. By similarity, we have $|N_G (x_3^3)\cap \{x_2^1,x_2^2\}|\geq 1$ and $|N_G (x_3^2)\cap \{x_3^1,x_1^3\}|\geq 1$, see Figure \ref{fi2}. In any case, we have $3K_ 2\subseteq G$; a contradiction again.
 	\begin{figure}[ht] 
 		\begin{tabular}{ccc}
 			\begin{tikzpicture}
 			\node [draw, circle, fill=black, inner sep=1.3pt, label=below:$x_3^2$] (x_3^2) at (0,0) {};
 			\node [draw, circle, fill=black, inner sep=1.3pt, label=below:$x_3^1$] (x_3^1) at (1,0) {};
 			\node [draw, circle, fill=black, inner sep=1.3pt, label=below:$x_2^2$] (x_2^2) at (2,0) {};
 			\node [draw, circle, fill=black, inner sep=1.3pt, label=below:$x_3^3$] (x_3^3) at (3,0) {};
 			
 			\node [draw, circle, fill=black, inner sep=1.3pt, label=left:$x_1^3$] (x_1^3) at (.5,0.85) {};
 			\node [draw, circle, fill=black, inner sep=1.3pt, label=right:$x_2^1$] (x_2^1) at (2.5,0.85) {};
 			
 			\node [draw, circle, fill=black, inner sep=1.3pt, label=left:$x_1^1$] (x_1^1) at (1,1.76) {};
 			\node [draw, circle, fill=black, inner sep=1.3pt, label=right:$x_1^2$] (x_1^2) at (2,1.76) {};

 			\node [draw, circle, fill=black, inner sep=1.3pt, label=above:$x_2^3$] (x_2^3) at (1.5,2.5) {};

 			\draw  (x_1^1)--(x_1^2)--(x_2^1)--(x_2^2)--(x_3^1)--(x_1^3)--(x_1^1)--cycle;
 			\draw [line width=1.pt,dash pattern=on 2pt off 2pt](x_1^1)--(x_2^3) (x_3^1)--(x_3^2) (x_2^1)--(x_3^3);
 			
 			\end{tikzpicture}\\
 		\end{tabular}\\
 		\caption{$B=(3,2,1)$}
 		\label{fi2}
 	\end{figure}

 	{\bf Case 2:} $B=(2,2,2)$.  W.l.g let $ Y_i=\{x_1^i,x_2^i\}$ for $i=1,2,3$. In this case, we have $C_6\cong w_1w_2w_3w_4w_5w_6w_1$. W.l.g let $w_1=x_1^1, w_2=x_1^2$. Since $|Y_3|=2$ and $w_4w_5\in E(C_6)$ we have $|\{w_3,w_6\}\cap Y_3|\geq 1$.  If $|\{w_3,w_6\}\cap Y_3|=2$, then considering  Figure \ref{fi3}(a), the proof is same  as case 1. So let $|\{w_3,w_6\}\cap Y_3|=1$.  W.l.g let $w_3=x_1^3$, $x_2^3=w_5$, $x_2^1=w_4 , x_2^2=w_6$. In this case,  consider  Figure \ref{fi3}(b) and the proof is  same as case 1. Hence, in any case, we have $3K_ 2\subseteq G$; again a contradiction.
 	
 	\begin{figure}[ht] 
 		\begin{tabular}{ccc}
 			\begin{tikzpicture}
 			\node [draw, circle, fill=black, inner sep=1.3pt, label=below:$x_3^2$] (x_3^2) at (0,0) {};
 			\node [draw, circle, fill=black, inner sep=1.3pt, label=below:$x_2^1$] (x_2^1) at (1,0) {};
 			\node [draw, circle, fill=black, inner sep=1.3pt, label=below:$x_2^2$] (x_2^2) at (2,0) {};
 			\node [draw, circle, fill=black, inner sep=1.3pt, label=below:$x_3^1$] (x_3^1) at (3,0) {};
 			
 			\node [draw, circle, fill=black, inner sep=1.3pt, label=left:$x_2^3$] (x_1^3) at (.5,0.85) {};
 			\node [draw, circle, fill=black, inner sep=1.3pt, label=right:$x_1^3$] (x_2^3) at (2.5,0.85) {};
 			
 			\node [draw, circle, fill=black, inner sep=1.3pt, label=left:$x_1^1$] (x_1^1) at (1,1.76) {};
 			\node [draw, circle, fill=black, inner sep=1.3pt, label=right:$x_1^2$] (x_1^2) at (2,1.76) {};
 			
 			\node [draw, circle, fill=black, inner sep=1.3pt, label=above:$x_3^3$] (x_3^3) at (1.5,2.5) {};
 			
 			\draw  (x_1^1)--(x_1^2)--(x_2^3)--(x_2^2)--(x_2^1)--(x_1^3)--(x_1^1)--cycle;
 			\draw (x_1^2)--(x_3^3) (x_3^1)--(x_2^2) (x_3^2)--(x_1^3);
 			\draw [line width=1.pt,dash pattern=on 2pt off 2pt](x_1^1)--(x_3^3) (x_3^1)--(x_2^3) (x_2^1)--(x_3^2);
 			\end{tikzpicture}
 			&&$~~~~~~~~~~~~~~~~~~~~$
 			\begin{tikzpicture}
 			\node [draw, circle, fill=black, inner sep=1.3pt, label=below:$x_3^1$] (x_3^1) at (0,0) {};
 			\node [draw, circle, fill=black, inner sep=1.3pt, label=below:$x_2^3$] (x_2^3) at (1,0) {};
 			\node [draw, circle, fill=black, inner sep=1.3pt, label=below:$x_2^1$] (x_2^1) at (2,0) {};
 			\node [draw, circle, fill=black, inner sep=1.3pt, label=below:$x_3^2$] (x_3^2) at (3,0) {};
 			
 			\node [draw, circle, fill=black, inner sep=1.3pt, label=left:$x_2^2$] (x_2^2) at (.5,0.85) {};
 			\node [draw, circle, fill=black, inner sep=1.3pt, label=right:$x_1^3$] (x_1^3) at (2.5,0.85) {};
 			
 			\node [draw, circle, fill=black, inner sep=1.3pt, label=left:$x_1^1$] (x_1^1) at (1,1.76) {};
 			\node [draw, circle, fill=black, inner sep=1.3pt, label=right:$x_1^2$] (x_1^2) at (2,1.76) {};
 			
 			\node [draw, circle, fill=black, inner sep=1.3pt, label=above:$x_3^3$] (x_3^3) at (1.5,2.5) {};
 			
 			\draw  (x_1^1)--(x_1^2)--(x_1^3)--(x_2^1)--(x_2^3)--(x_2^2)--(x_1^1)--cycle;
 			\draw (x_2^3)--(x_3^1) (x_2^1)--(x_3^2) (x_1^1)--(x_3^3);
 			\draw [line width=1.pt,dash pattern=on 2pt off 2pt](x_2^2)--(x_3^1) (x_1^3)--(x_3^2) (x_1^2)--(x_3^3);
 			\end{tikzpicture}\\
 			$a$&&$b$
 		\end{tabular}
 		\caption{(a) $|\{w_3,w_6\}\cap Y_3|=2$,$~~~~~~~$ (b) $|\{w_3,w_6\}\cap Y_3|=1$}
 		\label{fi3}
 	\end{figure}
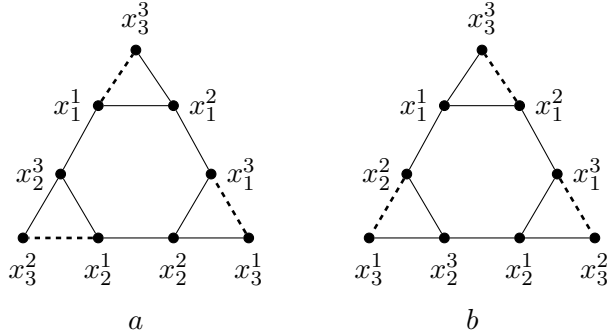
 	
 	By cases $1$ and $2$ we have $3K_ 2\subseteq G$. Thus the proof is complete and  the proposition  holds.
 \end{proof}
  We determined the exact value of the multipartite Ramsey number  $m_3(nK_2,C_7) $ for $n\geq 3$ in the following theorem:
 
 \begin{theorem}\label{vth1}
 	For each $n\geq 3$ we have $m_3(nK_2,C_7) = n$.
 \end{theorem}
 \begin{proof} 
 	First we show that $m_3(nK_2,C_7) \geq n$. Consider the coloring given by $K_{3\times (n-1)}=G^r\cup G^b$ where $G^r\cong K_{n-1, n-1}$ and $G^b\cong  K_{n-1,2(n-1)}$. Since $|V(G^r)|=2(n-1)$ and $G^b$ is  bipartite, we have $nK_2\not\subseteq G^r$ and $C_7\not\subseteq G^b$, that is, $m_3(nK_2,C_7) \geq n$. For the upper bound, consider  $K_{3\times n}$ with partite sets $X_i=\{x_1^i,x_2^i,\ldots, x_n^i\}$ for $i=1,2,3$. We will prove this by  induction. For $n=3$, by Proposition \ref{pr1}, the theorem holds. Suppose that $m_3(nK_2,C_7) \leq n$ for each $ n\geq 4$. We will show that $m_3((n+1)K_2,C_7) \leq n+1$ as follows: By contradiction we may assume that  $m_3((n+1)K_2,C_7) > n+1$, that is, $K_{3\times (n+1)}$ is $2$-colorable to $((n+1)K_2,C_7)$, say $(n+1)K_2\not\subseteq G$ and $C_7\not\subseteq\overline{G}$. Let $X'_i=X_i\setminus \{x^i_1\}$. Hence  by the induction hypothesis  we have $m_3(nK_2,C_7) \leq n$. Therefore since $|X'_i|=n$ and $C_7\not\subseteq \overline{G}[X'_1,X'_2,X'_3]$ we have $M=nK_2\subseteq G[X'_1,X'_2,X'_3]$. If there exists $i$ and $j$ such that $x_1^ix_1^j\in E(G)$, then we have $(n+1)K_2\subseteq G$; a contradiction. Hence we have  $x_1^ix_1^j\in E(\overline{G})$ for $i,j\in \{1,2,3\}$. Let $A=V(K_{3\times n}) \setminus V(M)$. Hence we have $|A|=3n-2n=n$. Since $(n+1)K_2\not\subseteq G$ we have $G[A,x^1_1,x^2_1,x^3_1]\subseteq \overline{G}$. Since $|A|=n\geq 4$ one can easily check that, in any case,  we have $H\subseteq \overline{G}$, where $H\in \{K_{5,1,1}, K_{4,2,1}, K_{3,3,1},K_{3,2,2}\}$. If  $H\in \{K_{3,3,1},K_{3,2,2}\}$, one can easily check that we have $C_7\subseteq H\subseteq \overline{G}$; a contradiction again. So let $H\in \{K_{5,1,1}, K_{4,2,1}\}$ and consider the following cases:
 	
 	\bigskip
 	{\bf Case 1:} $ A\subseteq X_i$  for only one $i$, that is, $H=K_{5,1,1}$.  W.l.g let $A\subseteq X_1$ and $\{x^1_2,x^1_3,\ldots,x^1_5\}\subseteq A$. Then we have $K_{n+1,1,1}\subseteq \overline{G}$ and $M\subseteq G[X_2,X_3]$. Since $n\geq 4$, we have $|M|\geq 4$, that is, there exists at least two edges, say $e_1= x_1y_1$ and $e_2=x_2y_2 $ in $E(M)$, where $\{x_1,x_2,y_1,y_2\} \subseteq X_2\cup X_3$. W.l.g let $|N_G(x_i)\cap A|\geq |N_G(y_i)\cap A|$ for $i=1,2$. One can easily check that $|N_G(y_i)\cap A|\leq 1$, otherwise we have $(n+1)K_2\subseteq G$; a contradiction. Since $|N_G(y_i)\cap A|\leq 1$ and $|A|\geq 5$, we have $|N_{\overline{G}}(y_i)\cap A|\geq 4$. Hence we have $|N_{\overline{G}}(y_1)\cap N_{\overline{G}}(y_2)\cap A|\geq 3$. W.l.g we may assume that $\{x_1^1, x_2^1, x_3^1\}\subseteq N_{\overline{G}}(y_1)\cap N_{\overline{G}}(y_2)\cap A$. In this case, we have $C_7\subseteq \overline{G}[x_1^1, x_2^1, x_3^1, x^2_1,x^3_1,y_1,y_2]\subseteq \overline{G}$; a contradiction again.
 	
 	\bigskip
 	{\bf Case 2:} $H= K_{4,2,1}$.  W.l.g let $|A\cap  X_1|=n-1$ and $|A\cap X_2|=2$. Let $\{x^1_2,x^1_3,\ldots,x^1_4\}\subseteq A\cap  X_1$ and $x_2^2\in  A\cap X_2$, that is, we have $K_{4,2,1}\subseteq K_{n,2,1}=G[A,x_1^1, x^1_2,x^1_3]\subseteq \overline{G}$ and $M\subseteq K_{1,n-1,n}$. That is, there exists at least one edge, say $e= xy$, where $x\in X_2$ and $y\in X_3$. W.l.g let $|N_G(x)\cap A|\geq |N_G(y)\cap A|$. One can easily check that $|N_G(y)\cap A|\leq 1$. Hence we have $|N_{\overline{G}}(y)\cap A|\geq3$ and the proof is same as case 1.
 	
 	\bigskip
 	By cases 1 and 2 we have the assumption that $m_3((n+1)K_2,C_7) > n+1$ does not hold. Now we have $m_3(nK_2,C_7) = n$ for each $n\geq 3$. This completes the induction step and the proof.
 \end{proof}
 
 \begin{lemma}\label{le3}
 	For  $j\geq 3$ and $n\geq j$ we have $m_j(nK_2,C_7) \geq \lceil \frac{2n+2}{j}\rceil$.
 \end{lemma}
 \begin{proof}
 	To show that  $m_j(nK_2,C_7) \geq \lceil \frac{2n+2}{j}\rceil$, assume that $ \lceil \frac{2n+2}{j}\rceil\geq 1$. Consider the coloring given by $K_{j\times t_0}=G^r\cup G^b$ where $t_0=\lceil \frac{2n+2}{j}\rceil -1$ such that $G^r\cong  K_{(j-1)\times t_0}$ and  $G^b\cong K_{t_0,(j-1)t_0}$. Since  $G^b$ is bipartite, we have $C_7\not\subseteq G^b$, and
 	\[|V(G^r)|=(j-1)\times t_0=  (j-1)(\lceil \frac{2n+2}{j}\rceil -1)=(j-1)(\lceil \frac{2n+2}{j}\rceil)-(j-1)\]
 	\[\leq (j-1)(\frac{2n+2}{j}+1)-(j-1)=j\times (\frac{2n+2}{j}) -\frac{2n+2}{j}.\]
 	Since $n\geq j$ we have $|V(G^r)|<2n$. Hence we have $nK_2\not\subseteq G^r$. Since $K_{j\times t_0}=G^r\cup G^b$ we have $m_j(nK_2,C_7) \geq \lceil \frac{2n+2}{j}\rceil$ for $n\geq j\geq 3$. 
 \end{proof}
 \begin{lemma}\label{le4}
 	$m_4(4K_2,C_7) = 3$.
 \end{lemma}
 \begin{proof}
 	By Lemma \ref{le3} we have $m_4(4K_2,C_7) \geq 3$. For the upper bound  consider the coloring given by $K_{4\times 3}=G^r\cup G^b$ such that $C_7\not\subseteq G^b$. Since $m_3(3K_2,C_7) = 3$, we have $3K_2\subseteq G^r[X_1,X_2,X_3]\subseteq G^r$. Let $M=3K_2$, hence we have $|V(X_1\cup X_2\cup X_3)\setminus V(M)|=3$. W.l.g let $A=\{w_1,w_2,w_3\}$ be this vertices. If $E(G^r)\cap E(G[X_4,A])\neq  \emptyset$ then we have $4K_2\subseteq G^r$. So let $ K_{3,3}\subseteq G[X_4,A]\subseteq G^b$. Consider the edge $e=v_1v_2 \in E(M)$,  it is easy to show that  $|N_{G^b}(v_i)\cap X_4|\geq 2$ for some $i\in\{1,2\}$, otherwise we have $4K_2\subseteq G^r$. In any case, one can easily check that  $C_7\subseteq G^b$; a contradiction. Thus giving $m_4(4K_2,C_7) = 3$.
 	
 \end{proof}
 
 \begin{lemma}\label{le5}
 	For $n\geq 4$ we have $m_4(nK_2,C_7)= \lceil \frac{n+1}{2}\rceil$.
 \end{lemma}
 \begin{proof}
 	By Lemma \ref{le3} we have  $m_4(nK_2,C_7) \geq \lceil \frac{n+1}{2}\rceil$. To prove $m_4(nK_2,C_7) \leq \lceil \frac{n+1}{2}\rceil$, consider  $K_{4\times t}$ with  partite set $X_i=\{x_1^i,x_2^i,\ldots, x_t^i\}$ for $i=1,2,3,4$, where $t=\lceil \frac{n+1}{2}\rceil$. We will prove this by  induction. For $n=4$ by Lemma  \ref{le4},  the lemma holds. Now we consider the following cases:
 	
 	\bigskip
 	{\bf Case 1:} $n=2k$, where $k\geq 3$. Suppose that $m_4(n'K_2,C_7) \leq \lceil \frac{n'+1}{2}\rceil $ for each $n'< n$. We will show that $m_4(nK_2,C_7) \leq \lceil \frac{n+1}{2}\rceil$ as follows: By contradiction, we may assume that  $m_4(nK_2,C_7) > \lceil \frac{n+1}{2}\rceil$, that is, $K_{4\times t}$ is $2$-colorable to $(nK_2,C_7)$, say $nK_2\not\subseteq G$ and $C_7\not\subseteq\overline{G}$. Let $X'_i=X_i\setminus \{x^i_1\}$ for $i=1,2,3,4$. Hence  by the induction hypothesis,  we have $m_4((n-1)K_2,C_7) \leq \lceil \frac{n}{2}\rceil=k$. Therefore, since $|X'_i|=k=\frac{n}{2}$ and $C_7\not\subseteq \overline{G}$ we have $M=(n-1)K_2\subseteq G[X_1',X_2',X_3',X_4']$. If there exists $i,j\in \{1,2,3,4\}$, where $x_1^ix_1^j\in E(G)$, then  $nK_2\subseteq G$; a contradiction. Now we have  $K_4 \cong \overline{G}[x_1^1,x_1^2,x_1^3,x_1^4]\subseteq \overline{G^g}$. Since  $nK_2\not\subseteq G$ and $\lceil \frac{n+1}{2}\rceil=\lceil\frac{2k+1}{2}\rceil=k+1$, we have $|V(K_{4\times k})\setminus V(M)|=2n-2(n-1)=2$, that is, there exists two vertices, say $w_1$ and $w_2$ in $V(K_{4\times k}) \setminus V(M)$. Since $nK_2\not\subseteq G$ we have $G[S]\subseteq \overline{G}$,  where $S=\{x_1^i~|~i=1,2,3,4\}\cup \{w_1,w_2\}$. Hence we have the following claim:
 	
 	\begin{claim}\label{c8}
 		Let $e=v_1v_2\in E(M)$ and w.l.g we may assume that $|N_{G}(v_1)\cap S|\geq |N_{G}(v_2)\cap S| $. If $|N_{G}(v_1) \cap S|\geq 2$ then $|N_{G}(v_2)\cap S|=0$. If $|N_{G}(v_1)\cap S|=1$ then $|N_{G}(v_2)\cap S|\leq 1$. If $|N_{G}(v_i)\cap S|=1$  then $v_1$ and $v_2$ have the same neighbour in $S$.
 	\end{claim}
{\it Proof.} By contradiction. We may assume that $\{w,w'\}\subseteq N_{G}(v_1)\cap S$ and $w''\in N_{G}(v_2)\cap S$, in this case, we set $ M' = (M\setminus \{v_1v_2\}) \cup\{v_1w,v_2w''\}$. Clearly, $M'$ is a matching whith $|M'| > |M|=n-1$, with contradicts the  $nK_2\not\subseteq G$ . If $|N_{G}(v_i)\cap S|=1$ and $v_i$ has a different neighbour then the proof is  same.
 
 \bigskip
 	Since $n\geq 4$ and $|M|\geq 3$. If $\{w_1,w_2\}\subseteq X_i$, say $ X_1$, then there is at least one edge, say $e=vu\in E(M)$ such that $v,u\notin X_1$. Otherwise, we have $C_7\subseteq K_{3\times 3} \subseteq \overline{G}[X_2,X_3,X_4]$; again a contradiction. W.l.g let $|N_{G}(v)\cap S|\geq |N_{G}(u)\cap S|$. Now  by Claim \ref{c8} we have $|N_{G}(u)\cap S| \leq 1$. One can easily check that in any case,  we have $C_7\subseteq \overline{G}[S\cup\{u\}]$; again a contradiction. So w.l.g let $w_1\in X_1$ and $w_2\in X_2$. In this case, since $|N_{G}(u)\cap S| \leq 1$, we have  $C_7\subseteq \overline{G}[S\cup\{u\}]$; a contradiction again.
 	
 	\bigskip
 	{\bf Case 2:} $n=2k+1$ where $k\geq 2$, $|X_i|=k+1$. Suppose that $m_4((n-2)K_2,C_7) \leq \lceil \frac{n-2+1}{2}\rceil $ for $n\geq 2$. We  show that $m_4(nK_2,C_7) \leq \lceil \frac{n+1}{2}\rceil$ as follows: By contradiction, we may assume that  $m_4(nK_2,C_7) > \lceil \frac{n+1}{2}\rceil$, that is, $K_{4\times t}$ is $2$-colorable to $(nK_2,C_7)$, say $nK_2\not\subseteq G$ and $C_7\not\subseteq\overline{G}$. Let $X'_i=X_i\setminus \{x^i_1\}$.  By the induction hypothesis,   we have $m_4((n-2)K_2,C_7) \leq \lceil \frac{n-1}{2}\rceil=\lceil \frac{2k}{2}\rceil=k$. Therefore, since $|X'_i|=k$ and $C_7\not\subseteq \overline{G}$ we have $M=(n-2)K_2\subseteq G[X'_1,X'_2,X'_3,X'_4]$ and thus we have the following claim:
 	
 	\begin{claim}\label{c9}
 		There exist two edges, say $e_1=uv $ and $e_2=u'v'$ in $E(M)=E((n-2)K_2)$, such that $v,v',u$ and $u'$  are in different partite.
 	\end{claim}
{\it Proof.} W.l.g assume that  $v\in X'_1$ and $u\in X'_2$. By contradiction, assume that $|E(M)\cap E(G[X'_3,X'_4])|=0$, that is, $G[X'_3,X'_4]\subseteq \overline{G}$. Since $|V(M)|=2(n-2)$ and $|X'_i|=k$, we have $|V(M)\cap X'_i|\geq k-2$. Since $k\geq 3$, $|V(M)\cap X'_j|\geq 1~(j=3,4)$. W.l.g let $e_j'=x_jy_j\in E(M)$ where $x_j\in V(M)\cap X'_j$. And w.l.g we may assume that $y_3\in V(M)\cap X'_1$. Hence we have $y_4\in V(M)\cap X'_1$.  In other words, take $e_1=x_3y_3$ and $e_2=x_4y_4$ and the proof is complete. Hence we have $|E(M)\cap E(G[X'_2,X'_j])|=0$ for $j=3,4$, in other words, if there exists $e''\in E(M)\cap E(G[X'_2,X'_j])$, then set $e_1=e_1'$ and $e_2=e''$ and the proof is complete. Therefore, for each $e\in E(M)$ we have $v(e)\cap X_1'\neq \emptyset$ which means that  $|M|\leq X'_1=k$; a contradiction to $|M|$.
  
  \bigskip
 	Now by Claim \ref{c9} there exist two edges, say $e_1=uv $ and $e_2=u'v'$ in $E(M)=E((n-2)K_2)$, such that $v,v',u$ and $u'$  are in  different partite.  W.l.g let $e_1=x_1x_2$ and $e_2=x_3x_4$ be this edges and let $x_i\in X_i'$ for $i=1,2,3,4$. Set $X''_i=X_i\setminus\{x_i\}$, hence we have $|X''_i|=k$. Since $C_7\not\subseteq \overline{G}$ we have $C_7\not\subseteq \overline{G}[X_1'',X_2'',X_3'',X''_4]$. Therefore,  by the induction hypothesis, we have $(n-2)K_2\subseteq G[X_1'',X_2'',X_3'',X''_4]$. Let $M=(n-2)K_2\subseteq G[X_1'',X_2'',X_3'',X''_4]$,   set $M^*=M\cup \{e_1,e_2\}$ hence $|M^*|=n$, that is, $nK_2\subseteq G$; again  a contradiction.  Hence the assumption that $m_4(nK_2,C_7) > \lceil \frac{n+1}{2}\rceil$ dose not hold and we have $m_4(nK_2,C_7) \leq \lceil \frac{n+1}{2}\rceil$. This completes the induction step and the proof is complete.
 	By cases 1, 2 we have $m_4(nK_2,C_7) = \lceil \frac{n+1}{2}\rceil$ for $n\geq 4$.
 \end{proof}
  The results of  Proposition \ref{pr1}, Lemmas  \ref{le4} and \ref{le5} and  Theorem  \ref{vth1} concludes the proof of  Theorem \ref{th2}.
% \begin{theorem}
 %	Let $j\in\{ 2,3,4\}$ and $n\geq 2$. Then 
%	\[
%	m_j(nK_2,C_7)= \left\lbrace
%	\begin{array}{ll}
%	\infty & ~~~~~~j=2,n\geq2 ,~\vspace{.2 cm}\\
%	2 & ~~~~~~(j,n)=(4,2),~\vspace{.2 cm}\\
%	3 & ~~~~~~(j,n)=(3,2), (4,3), ~\vspace{.2 cm}\\
%	n & ~~~~~~j=3,n\geq 3 ,~\vspace{.2 cm}\\
%	\lceil \frac{n+1}{2} \rceil &  ~~~~~~ j=4,n\geq 4 .
%	\end{array}
%	\right.
%	\]
 %\end{theorem}
%%%%%%%%%%%%%%%%%%%%%%%%%%%%%%%%%%%%%%%%%
\bibliographystyle{spmpsci} 
\bibliography{yas-Corrected}

\begin{thebibliography}{1}
\providecommand{\url}[1]{{#1}}
\providecommand{\urlprefix}{URL }
\expandafter\ifx\csname urlstyle\endcsname\relax
  \providecommand{\doi}[1]{DOI~\discretionary{}{}{}#1}\else
  \providecommand{\doi}{DOI~\discretionary{}{}{}\begingroup
  \urlstyle{rm}\Url}\fi

\bibitem{bondy1976graph}
Bondy, J.A., Murty, U.S.R.: Graph theory with applications.
\newblock American Elsevier Publishing Co., Inc., New York (1976)

\bibitem{burger2004diagonal}
Burger, A.P., Grobler, P.J.P., Stipp, E.H., van Vuuren, J.H.: Diagonal {R}amsey
  numbers in multipartite graphs.
\newblock Util. Math. \textbf{66}, 137--163 (2004)

\bibitem{burger2004ramsey}
Burger, A.P., Vuuren, J.H.V.: Ramsey numbers in complete balanced multipartite
  graphs. {II}. {S}ize numbers.
\newblock Discrete Math. \textbf{283}(1-3), 45--49 (2004).
\newblock \doi{10.1016/j.disc.2004.02.003}.
\newblock
  \urlprefix\url{https://doi.org.creativaplus.uaslp.mx/10.1016/j.disc.2004.02.003}

\bibitem{erdos2009partition}
Erd\"{o}s, P., Rado, R.: A partition calculus in set theory.
\newblock Bull. Amer. Math. Soc. \textbf{62}, 427--489 (1956).
\newblock \doi{10.1090/S0002-9904-1956-10036-0}.
\newblock
  \urlprefix\url{https://doi.org.creativaplus.uaslp.mx/10.1090/S0002-9904-1956-10036-0}

\bibitem{jayawardene2016size}
Jayawardene, C., Baskoro, E.T., Samarasekara, L., Sy, S.: Size multipartite
  {R}amsey numbers for stripes versus small cycles.
\newblock Electron. J. Graph Theory Appl. (EJGTA) \textbf{4}(2), 157--170
  (2016).
\newblock \doi{10.5614/ejgta.2016.4.2.4}.
\newblock
  \urlprefix\url{https://doi.org.creativaplus.uaslp.mx/10.5614/ejgta.2016.4.2.4}

\bibitem{luczak2018multipartite}
\L.~uczak, T., Polcyn, J.: The multipartite {R}amsey number for the 3-path of
  length three.
\newblock Discrete Math. \textbf{341}(5), 1270--1274 (2018).
\newblock \doi{10.1016/j.disc.2018.01.015}.
\newblock
  \urlprefix\url{https://doi.org.creativaplus.uaslp.mx/10.1016/j.disc.2018.01.015}

\bibitem{lusiani2017size}
Lusiani, A., Baskoro, E.T., Saputro, S.W.: On size multipartite {R}amsey
  numbers for stars versus paths and cycles.
\newblock Electron. J. Graph Theory Appl. (EJGTA) \textbf{5}(1), 43--50 (2017).
\newblock \doi{10.5614/ejgta.2017.5.1.5}.
\newblock
  \urlprefix\url{https://doi.org.creativaplus.uaslp.mx/10.5614/ejgta.2017.5.1.5}

\bibitem{yek}
Perondi, P.H., Carmelo, E.L.M.: Set and size multipartite {R}amsey numbers for
  stars.
\newblock Discrete Appl. Math. \textbf{250}, 368--372 (2018).
\newblock \doi{10.1016/j.dam.2018.05.016}.
\newblock
  \urlprefix\url{https://doi.org.creativaplus.uaslp.mx/10.1016/j.dam.2018.05.016}

\bibitem{sy2011size}
Sy, S.: On the size multipartite {R}amsey numbers for small path versus
  cocktail party graphs.
\newblock Far East J. Appl. Math. \textbf{55}(1), 53--60 (2011)

\end{thebibliography}
\end{document}